\documentclass[12pt,a4paper]{article}

\usepackage{amsmath}
\usepackage{amsthm}

\usepackage[latin1]{inputenc}
\usepackage[T1]{fontenc}

\usepackage[scaled=0.92]{helvet}
\usepackage{courier}

\usepackage{eufrak}

\newtheorem{theorem}{Theorem}

\newtheorem{lemma}[theorem]{Lemma}
\newtheorem{proposition}[theorem]{Proposition}

\theoremstyle{definition}

\theoremstyle{remark}

\DeclareMathOperator{\rank}{rank}

\begin{document}

\title{Extending $\pi$-systems to bases of root systems}

\author{Helmer ASLAKSEN\\
Department of Mathematics\\
National University of Singapore\\
Singapore 117543\\
Singapore\\[6pt]
{aslaksen@math.nus.edu.sg}\\
{www.math.nus.edu.sg/aslaksen/}
\and
Mong Lung LANG\\
Department of Mathematics\\
National University of Singapore\\
Singapore 117543\\
Singapore\\[6pt]
{matlml@math.nus.edu.sg}}


\date{}

\maketitle

\begin{abstract}
Let $R$ be an indecomposable root system. It is well known that any
root is part of a basis $B$ of $R$. But when can you extend a set of
two or more roots to a basis $B$ of $R$? A $\pi$-system is a
linearly independent set of roots, $C$, such that if $\alpha$ and
$\beta$ are in $C$, then $\alpha - \beta$ is not a root. We will use
results of Dynkin and Bourbaki to show that with two exceptions,
$A_3 \subset B_n$ and $A_7 \subset E_8$, an indecomposable
$\pi$-system whose Dynkin diagram is a subdiagram of the Dynkin
diagram of $R$ can always be extended to a basis of $R$.
\end{abstract}

\section{Introduction}

Let $R$ be an indecomposable root system in a Euclidean space $V$. A
subset $B$ of $R$ is called a basis of $R$ if $B$ is a vector space
basis of $V$ and each root of $R$ can be written as a linear
combination of roots in $B$ with integral coefficients that are all
nonnegative or all nonpositive. It is well known that any root is
part of a basis $B$ of $R$. But when can you extend a set of two or
more roots to a basis $B$ of $R$? A \emph{$\pi$-system}
(\cite{Dyn,Oni-Vin}) is a linearly independent set of roots, $C$,
such that if $\alpha$ and $\beta$ are in $C$, then $\alpha - \beta$
is not a root. (It is not assumed to be linearly independent in
\cite{Oni-Vin}.) A subset of a basis will be a $\pi$-system, and a
$\pi$-system will be a basis of a root subsystem. We can associate a
Dynkin diagram to a $\pi$-system, and in order to extend to a basis
of $R$, the Dynkin diagram of the $\pi$-system must be a subdiagram
of the Dynkin diagram of $R$. By a subdiagram we mean a diagram
obtained by deleting some nodes and their corresponding links. We
will assume that the nodes corresponding to short roots are marked,
and that the subdiagram preserves the marking. Hence two orthogonal
short roots do not form a subdiagram of $B_n$, while two orthogonal
long roots do. We will use results of Dynkin (\cite{Dyn} and
Bourbaki (\cite{Bou} to show that with two exceptions, $A_3 \subset
B_n$ and $A_7 \subset E_8$, an indecomposable $\pi$-system whose
Dynkin diagram is a subdiagram of the Dynkin diagram of $R$ can
always be extended to a basis of $R$. Our techniques can easily
handle the decomposable case, too, but the results become more
tedious to state, and we feel that it would distract from the main
ideas of the paper.

We would like to thank professor Robert V. Moody for pointing out
the results from Bourbaki to us during his visit to Singapore.

\section{Results from Dynkin and Bourbaki}

If $C$ is a set of roots, then $[C]$ denotes the set of all roots in
$R$ that are linear combinations of the roots in $C$ with integer
coefficients. Let $\Pi$ and $\Pi'$ be $\pi$-systems. We will say
that $\Pi'$ is obtained from $\Pi$ by an \emph{elementary
transformation} if $\Pi'$ is obtained by adjoining the lowest root
to one of the indecomposable components of $\Pi$ and then removing
one root from that component.

We will first state three results due to Dynkin (\cite[Theo.\ 5.1,
5.2 and 5.3]{Dyn}).

\begin{proposition}
\label{Dyn1}

Let $C$ be a $\pi$-system in a root system, $R$. Then $[C]$ is a
root subsystem of $R$ with basis $C$.
\end{proposition}

\begin{proposition}
\label{Dyn2}

Let $C$ be a $\pi$-system in an indecomposable root system, $R$, of
rank $n$. Then $C$ can be extended to a $\pi$-system, $D$, with $n$
elements.
\end{proposition}

\begin{proposition}
\label{Dyn3}

Let $D$ be a $\pi$-system with $n$ elements in an indecomposable
root system, $R$, of rank $n$. Then $D$ can be obtained by a
sequence of elementary transformations of a basis of $R$.
\end{proposition}

This shows that extending $C$ will not always give us a basis of
$R$, but that the extension can be obtained by a sequence of
elementary transformations of a basis of $R$.

The next three propositions are due to Bourbaki (\cite[Coro.\ to
Prop.\ 4, Coro.\ 4 to Prop.\ 20 and Prop.\ 24 of Chap.\ VI, \S
1]{Bou}).

\begin{proposition}
\label{Bou1}

Let $V'$ be a subspace of $V$ and let $V''$ be the subspace spanned
by $R' = R \cap V'$. Then $R'$ is a root system in $V''$.

\end{proposition}

Notice that $V'$ may intersect $R$ so that $\dim V'' < \dim V'$.

\begin{proposition}
\label{Bou2}

If $C$ a subset of a basis $B$ of $R$ and $V'$ the subspace of $V$
spanned by $C$, then $C$ is a basis of the root system $R' = R \cap
V'$.

\end{proposition}

Notice that $[C]$ uses integer coefficients, while $R'= R \cap V'$
includes rational coefficients. This proposition says that a if $C$
is a subset of a basis, then $[C] = R'$.

\begin{proposition}
\label{Bou3}

Let $C$ be a basis of $R' = R \cap V'$, where $V'$ is a subspace of
$V$. Then $C$ can be extended to a basis $B$ of $R$ and $R'$ is the
set of roots in $R$ that are linear combinations of elements of $C$.

\end{proposition}

This says that if $[C] = R'$, then $C$ can be extended to a basis of
$R$.

Then next result is a simple combination of Propositions~\ref{Dyn1}
to~\ref{Bou3}.

\begin{theorem}
\label{extend}

Let $C$ be a $\pi$-system in $R$ and let $V'$ be the subspace of $V$
spanned by $C$. Then $[C]$ is a root subsystem of $R' = R \cap V'$,
$C$ can be obtained from a basis of $R'$ by a sequence of elementary
transformations, and $C$ can be extended to a basis $B$ of $R$ if
and only if $C$ is a basis of $R'$, i.e., $[C] = R'$.

\end{theorem}

\begin{proof}
It follows from Proposition~\ref{Dyn1} that $[C]$ is a root
subsystem of $R' = R \cap V'$ with basis $C$. In general $[C] \not=
R'$, so we can only use Proposition~\ref{Dyn3} to say that $C$ is
obtained from a basis of $R'$ by a sequence of elementary
transformations. If $[C] = R'$, then it follows from
Proposition~\ref{Bou3} that $C$ can be extended to a basis of $R$.
The converse follows from Proposition~\ref{Bou2}.
\end{proof}

\section{Extension results}

It follows from Theorem~\ref{extend} that $C$ will extend to a basis
of $R$ unless there is a root system $R'$ such that
\begin{equation*}
[C] \subset R' \subset R,
\end{equation*}
where $\rank \, [C] = \rank R' < \rank R$, $[C] \not = R'$ and $C$
can be obtained from a basis of $R'$ by a sequence of elementary
transformations. We are for simplicity assuming that $[C]$ is
indecomposable, and hence $R'$ is also indecomposable. The next
lemma is proved by inspection of the extended Dynkin diagrams.

\begin{lemma}
\label{transf}

The only indecomposable roots systems $[C] \subset R'$ where $C$ can
be obtained by a sequence of elementary transformations of a basis
of $R'$ are listed below.

\begin{gather*}
 A_3 \subset B_3, \quad
 D_n \subset B_n \quad \mbox{for $n \geq 4$,} \quad
 A_7 \subset E_7, \\
 A_8 \subset E_8, \quad
 D_8 \subset E_8, \quad
 B_4 \subset F_4, \quad
 D_4 \subset F_4, \quad
 A_2 \subset G_2.
\end{gather*}

$D_4 \subset F_4$ is obtained by using two elementary
transformations, while all the other only require one.

\end{lemma}

Since our pairs must sit inside a root system, $R$, of higher rank,
and the diagram of $C$ must be a subdiagram of the diagram of $R$,
there are only two cases that satisfies our conditions, $A_3 \subset
B_3 \subset B_n$ and $A_7 \subset E_7 \subset E_8$.

The only way $C$ can fail to extend is if $C$ is obtained from an
$R'$ diagram by an elementary transformation so that $C$ includes
the lowest root of $R'$. We will see in the next theorem that in
that case we can use the fact that the lowest root is a linear
combination of all the roots in the basis to show that $[C] \not =
R'$.

\begin{theorem}
\label{theorem}

An indecomposable $\pi$-system, $C$, in $R$ can be extended to a
basis $B$ of $R$, unless $[C] \subset R$ is on of the following two
cases.

\begin{enumerate}

\item $A_3 \subset B_n$ for $n \geq 4$ and $C = \{ r_1, r_2, r_3 \}$
has Dynkin diagram

\begin{center}
 \begin{picture}(60,10)(0,0)
 \multiput(0,0)(30,0){3}{\circle{4}}
 \multiput(2,0)(30,0){2}{\line(1,0){26}}
 \put(-3,4){$r_1$}
 \put(27,4){$r_2$}
 \put(57,4){$r_3$}
 \end{picture}
\end{center}
and $(r_1 + 2 r_2 + r_3)/2$ is a root in $B_n$.

\item $A_7 \subset E_8$ and $C = \{ r_1, \dots, r_7 \}$ has Dynkin
diagram

\begin{center}
 \begin{picture}(180,15)
 \multiput(0,0)(30,0){7}{\circle{4}}
 \multiput(2,0)(30,0){6}{\line(1,0){26}}
 \put(-3,4){$r_1$}
 \put(27,4){$r_2$}
 \put(57,4){$r_3$}
 \put(87,4){$r_4$}
 \put(117,4){$r_5$}
 \put(147,4){$r_6$}
 \put(177,4){$r_7$}
 \end{picture}
\end{center}
and
\[
(r_1 + 2 r_2 + 3 r_3 + 4 r_4 + 3 r_5 + 2 r_6 + r_7)/2
\]
is a root in $E_8$.

\end{enumerate}

\end{theorem}

\begin{proof}

We cannot use the pairs in Theorem~\ref{transf} involving $E_8$,
$F_4$ and $G_2$, because we cannot fit them into any bigger root
systems. We also cannot fit a $D_k$ diagram inside $B_n$, so we are
left with $A_3 \subset B_3 \subset B_n$ and $A_7 \subset E_7 \subset
E_8$.

We will use the bases for the root systems listed in \cite{Bou} and
denote the lowest root by $\alpha_0$.

For $A_3 \subset B_n$, either $A_3$ is the Y-branch at the end of
the extended diagram of $B_3$ or $A_3$ is part of the diagram of
$B_n$. In the first case, $r_2 = \alpha_2$ and either $r_1$ or $r_3$
must be the lowest root $-(\alpha_1 + 2 \alpha_2 + 2 \alpha_3)$ and
the other must be $\alpha_1$. In either case, we can recover the
deleted root from the lowest root, giving us $(r_1 + 2 r_2 + r_3)/2
= - \alpha_3 = - e_3$, which is a root of $B_n$, but not in $[C]$.
So $V'$ is the span of $\{e_1, e_2, e_3\}$, and $R' = B_3$ while
$[C]=A_3$. Hence $C$ cannot be extended to a basis of $B_n$ by
Theorem~\ref{extend}.

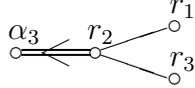
\begin{center}
 \begin{picture}(60,20)(0,-5)
 \multiput(0,0)(30,0){2}{\circle{4}}
 \put(2,-1){\line(1,0){26}}
 \put(2,1){\line(1,0){26}}
 \put(10,0){\line(2,1){10}}
 \put(10,0){\line(2,-1){10}}
 \put(-3,5){$\alpha_3$}
 \put(27,5){$r_2$}
 \put(32,1){\line(3,1){25}}
 \put(32,-1){\line(3,-1){25}}
 \put(60,-10){\circle{4}}
 \put(60,10){\circle{4}}
 \put(58,15){$r_1$}
 \put(58,-5){$r_3$}
 \end{picture}
\end{center}

However, if $C = \{ \alpha_1, \alpha_2, \alpha_3 \}$ corresponds to
an $A_3$ that is a subdiagram of the diagram of $B_n$, then $(r_1 +
2 r_2 + r_3)/2$ is not a root of $B_n$. In this case $V'$ does not
contain any short roots, so $R' = [C]$ and $C$ can be extended to a
basis of $B_n$.

For $A_7 \subset E_8$, the $A_7$ is either part of the diagram of
$E_8$ or is part of the extended diagram of $E_7$. In the $E_8$ case
we have

\[
 C = \{ \alpha_1, \dots, \alpha_8 \} - \{ \alpha_2 \},
\]

and in the extended $E_7$ case we have

\[
 C = \{ \alpha_0, \alpha_1, \dots \alpha_7 \} - \{ \alpha_2 \},
\]

where $\alpha_0 = - 2 \alpha_1 -2 \alpha_2 - 3 \alpha_3 - 4 \alpha_4
- 3 \alpha_5 - 2 \alpha_6 - \alpha_7$ is the lowest root of $E_7$.

In the extended $E_7$ case we don't know whether $r_1$ or $r_7$ is
the extended root, but again we can recover the deleted root from
the lowest root, giving us

\[
(r_1 + 2 r_2 + 3 r_3 + 4 r_4 + 3 r_5 + 2 r_6 + r_7)/2 = - \alpha_2.
\]

In the $E_8$ case this expression will not be a root. It follows
that in the extended $E_7$ case $C$ can be extended to a basis of
$E_8$, while in the $E_8$ case $C$ can only be extended to a basis
of $A_8$.
\end{proof}

Notice that $C = \{ \alpha_1, \alpha_2, \alpha_0, \} = \{ e_1-e_2,
e_2-e_3, -e_1-e_2 \}$, where $\alpha_0$ is the lowest short root in
$C_n$ is \emph{not} a $\pi$-system in $C_n$, since $e_1-e_2 -
(-e_1-e_2) = 2e_1$ is a root in $C_n$. We initially considered
linearly independent sets of roots with nonpositive inner products,
i.e., linearly independent, admissible (nonpositive inner product)
sets of roots instead of $\pi$-systems. $\pi$-systems are always
admissible, and for simply laced root systems, indecomposable
$\pi$-systems are admissible, since the only way the difference
between two roots can have the same length as the two roots is if
the angle between them is $\pi/3$. However, $C \subset C_n$ shows
that this is false for multiply laced root systems. In particular,
Exercise~34 on page~177 of \cite{Oni-Vin} appears to be incorrect.
(They do not require $\pi$-systems to be linearly independent, but
that does not make any difference.)

Notice that $D_k \subset B_k$ is listed in Lemma~\ref{transf}, while
$D_k \subset C_k$ is not. There are two standard ways of
constructing equal rank inclusions of root systems. One is to use
elementary transformations and is used by Borel and de Siebenthal
(\cite{Bor-Sib}). The other is to consider the set of short and long
roots in multiply laced root systems. $D_k$ forms the long roots in
$B_k$ and the short roots in $C_k$, but only the $B_k$ inclusion can
be obtained by an elementary transformation.

Notice also that we only talk about root systems, and not about Lie
subalgebras. Unless we know something about the Cartan subalgebras,
inclusions of root systems and inclusions of Lie algebras will not
necessarily correspond. The fact that $D_k \subset C_k$ does not
imply that $\mathfrak{so}(2n) \subset \mathfrak{sp}(2n)$.


\begin{thebibliography}{9}

\bibitem[Bor-Sib]{Bor-Sib}A. Borel and J. de Siebenthal, Les
sous-groupes fermés de rang maximum des groupes de Lie clos,
\emph{Comment.\ Math.\ Helvetici} \textbf{23} (1949), 200--221.

\bibitem[Bou]{Bou}N. Bourbaki, \emph{Groupes et algèbres de Lie,
Chapitres 4, 5 et 6, Éléments de mathématique}, Hermann, 1968.

\bibitem[Dyn]{Dyn}E.B. Dynkin, Semisimple Subalgebras of Semisimple
Lie Algebras, \emph{Mat.\ Sbornik N.S.} \textbf{30 (72)} (1952),
349--462 (Russian). English translation in American Mathematical
Society Translations, Series 2, vol.\ 6, \emph{Five Papers on
Algebra and Group Theory}, Amer.\ Math.\ Soc., 1957, pp.\ 111--245.

\bibitem[Oni-Vin]{Oni-Vin}A.L. Onishchik and E.B. Vinberg, \emph{Lie
Groups and Algebraic Groups}, Springer-Verlag, 1990.

\end{thebibliography}
\end{document}